\documentclass[11pt]{article} 
 
\usepackage[utf8]{inputenc} 

 \usepackage{array} 

\usepackage{mathtools}
\usepackage{amsfonts,hyperref,url}
\usepackage{amssymb}
\usepackage{amsthm}

\newcommand{\ba}{\begin{eqnarray}}
\newcommand{\ea}{\end{eqnarray}}

\newcommand{\Rr}{\mathbb{R}}
\newtheorem{theorem}{Theorem}[section]
\newtheorem{hyp}{Assumption}
\newtheorem{property}{Property}
\newtheorem{lemma}[theorem]{Lemma}
\newtheorem{proposition}[theorem]{Proposition}

\newtheorem{definition}[theorem]{Definition\rm}
\newtheorem{remark}{\it Remark\/}

\everymath{\displaystyle}

\begin{document}
 
\title{Pathwise integration with respect to paths of finite quadratic variation}
\author{Anna ANANOVA and Rama CONT}

\date{
Revised: August 2016. To appear in:  {\it Journal de Math\'ematiques Pures et Appliqu\'ees}.} 

\maketitle
\begin{abstract}
We study  a pathwise integral with respect to paths of finite quadratic variation, defined as the limit of non-anticipative Riemann sums for  gradient-type integrands.
We show that the integral  satisfies  a pathwise isometry property, analogous to the well-known Ito isometry for stochastic integrals. This property is then used to represent the integral as  a continuous map on an appropriately defined vector space of integrands.
Finally, we obtain a  pathwise 'signal plus noise' decomposition for regular functionals of an irregular path with non-vanishing quadratic variation, as a unique sum of a pathwise integral and a component with zero quadratic variation.
\end{abstract}
A. Ananova \& R. Cont (2016)  {\rm Pathwise integration with respect to paths of finite quadratic variation}, {\it Journal de Math\'ematiques Pures et Appliqu\'ees}.
\tableofcontents
\newpage

In his seminal paper 'Calcul d'Ito sans probabilit\'es' \cite{follmer1981}, Hans F\"ollmer  proved a change of variable formula for smooth functions of paths with infinite variation, using the concept of {\it quadratic variation along a sequence of partitions}.
A  path $\omega\in C^0([0,T],\mathbb{R})$ is said to have finite
quadratic variation along the sequence of partitions
$\pi^m=(0=t^{m}_0<t^m_1<\cdots<t^{m}_{k(m)}=T)$  if for any $t\in [0,T]$, the limit
\ba [\omega]_{\pi}(t):=\mathop{\lim}_{m\to\infty}\sum_{t^m_{i}\leq t} (\omega(t^m_{i+1})-\omega(t^m_i))^2 <\infty\nonumber\ea
exists and defines a continuous increasing function $[\omega]_\pi: [0,T]\to \mathbb{R}_+$ called the quadratic variation of $\omega$ along $\pi$. Extending this definition to vector-valued paths (see Section \ref{sec:FIC})
 F\"ollmer \cite{follmer1981} showed that, for integrands of the form  $\nabla f( \omega(t)), f\in C^2(\mathbb{R}^d)$  one may define a pathwise integral $\int \nabla f(\omega(t))d\omega(t)$  as a pointwise limit of  Riemann sums  along the sequence of partitions $\pi^m$.

This construction has been extended in various directions,  to less regular functions
\cite{bertoin1987,davis2014,perkowski2015b} and path-dependent functionals \cite{CF10B,cont2016,contriga2014}. In particular, Cont \& Fourni\'e \cite{CF10B} constructed pathwise integrals of the type
$$\int_0^T \nabla_\omega F(t,\omega)d^{\pi}\omega $$
  where $\nabla_\omega F$ is a directional derivative (Dupire derivative) of a non-anticipative functional $F$, and proved a change of variable formulas for such integrals.

This paper contributes several new results to the study of the pathwise approach to stochastic integration \cite{follmer1981} and its extension to path-dependent functionals \cite{CF10B}. 

Our main result  is an  isometry formula for the pathwise integral: we give conditions on the integrand $\phi $ and the path $\omega\in C^0([0,T],\mathbb{R})$   under which
$$  [\int_0^. \phi d^\pi\omega ]_\pi (t)  = \int_0^t |\phi|^2 d[\omega]_\pi. $$
Theorem \ref{theorem.isometry} gives a precise statement of this property for multidimensional paths.
Our conditions notably allow for continuous paths with H\"older exponent strictly less than $1/2$ and apply to  typical paths of Brownian motion and diffusion processes.

This result may be understood as a pathwise version of the well-known Ito isometry formula \cite{ikedawatanabe}, showing that the Ito isometry formula holds not just as an equality in expectation but in fact path by path.

This result  has several interesting consequences. 
First, it implies that the pathwise integral may be defined as a  a continuous mapping for the 'quadratic variation' metric on an appropriate space of integrands (Proposition \ref{prop.quotient}). This continuity property, together with its interpretation as a limit of Riemann sums, is what distinguishes our construction from other   pathwise constructions of stochastic integrals \cite{karandikar1995,nutz2012,perkowski2015b,russovallois} which  lack either a continuity property or an interpretation as a limit of non-anticipative Riemann sums.

A second consequence of the pathwise isometry property is a deterministic  'signal plus noise' decomposition for functionals of an irregular path: we show that any regular functional of path  $\bar\omega$ with non-vanishing quadratic variation may be uniquely decomposed as the sum of a pathwise integral with respect to $\bar\omega$ and a 'smooth' component  with zero quadratic variation (Proposition \ref{prop.decomposition}). This may be seen as a pathwise analogue of the  well-known semimartingale decomposition \cite{dm} (or, more precisely the decomposition of Dirichlet processes \cite{follmer1981b}). A  similar result was obtained by Hairer and Pillai \cite{hairerpillai2013}  in the context of rough path theory, using different techniques.

Finally, we clarify the pathwise nature of the integral defined in \cite{CF10B}: we show (Theorem \ref{theorem.isomuniq}) that this integral is indeed a pathwise limit of non-anticipative Riemann sums, which is important for interpretation and in applications.

{\bf Outline}  Section \ref{sec:FIC} recalls some key definitions and results on functional calculus from \cite{CF10B,cont2012} and recalls the definition of the F\"ollmer integral \cite{follmer1981} and its extension to path-dependent integrands by Cont \& Fourni\'e \cite{CF10B}. The 'isometry formula' for this integral is derived in Section \ref{sec:isometry} (Theorems \ref{theorem.isometry}).  In Section \ref{sec:LebesguePartition} we discuss the 'isometry formula' for Lebesgue partitions. We use the of Section \ref{sec:isometry} in Section \ref{sec:quotient} to represent the integral as a continuous map (Proposition \ref{prop.quotient}).  In Section \ref{sec:uniqueness} we clarify the 'pathwise' nature of this integral (Theorem \ref{theorem.isomuniq}) .
 Section \ref{sec:decomposition}   derives a  pathwise 'signal plus noise' decomposition for functionals of irregular paths (Proposition \ref{prop.decomposition}), which may be viewed as a deterministic analogue of the semimartingale decomposition for stochastic processes.

\section{Non-anticipative functional calculus and pathwise integration of gradient functionals}\label{sec:FIC}

Our approach relies on the Non-anticipative Functional Calculus \cite{CF10B,cont2012},  a functional calculus which applies to non-anticipative functionals of c\`adl\`ag paths with finite quadratic variation, in the sense of F\"ollmer \cite{follmer1981}.  We recall here some key concepts and results of this approach, following \cite{cont2012}.

 Let $ X $ be the canonical process on $ \Omega = D ( [0,T], \mathbb{R}^d ) $, and $ ( \mathcal{F}_t^0 )_{t \in [0, T]} $ be the filtration generated by $ X $. We are interested in {\it causal}, or {\it non-anticipative functionals} of $X$, that is, functionals  $ F:[0, T] \times D( [0, T], \mathbb{R}^d ) \mapsto \mathbb{R} $ such that
\begin{equation}  \forall \omega \in \Omega,\qquad F (t, \omega) = F (t, \omega_t),\label{eq:nonanticipative}\end{equation}
where $\omega_t:=\omega(t\wedge \cdot)$ is the stopped path of $\omega$ at time $t.$
The process $ t \mapsto F (t,X_t) $  then only depends on the path of $X$ up to $ t $ and is $ ( \mathcal{F}_t^0 ) $-adapted.

It is convenient to define such functionals  on the space of stopped paths \cite{cont2012}: a stopped path is an equivalence class in $ [0, T] \times D( [0, T], \mathbb{R}^d ) $ for the following equivalence relation:
\begin{equation} \label{equiv}
(t, \omega) \sim (t', \omega') \Longleftrightarrow \left( t=t' \quad \mathrm{and} \quad \omega(t\wedge .) = \omega'({t}\wedge  .)\quad \right).
\end{equation}
The space of stopped paths is defined as the quotient of $ [0, T] \times D( [0, T], \mathbb{R}^d ) $ by the equivalence relation \eqref{equiv}:
\[ \Lambda_T = \{ (t, \omega (t \wedge \cdot)), (t, \omega) \in [0, T] \times D( [0, T],\mathbb{R}^d ) \} = \left( [0, T] \times D( [0, T], \mathbb{R}^d ) \right) / \sim \]
We denote ${\mathcal W}_T$ the subset of $\Lambda_T$ consisting of continuous stopped paths.
We endow the set $\Lambda_T$ with a metric space structure by defining the following distance:
\begin{eqnarray*}
d_{\infty} ((t, \omega), (t', \omega')) = \sup_{u \in [0, T]} | \omega (u \wedge t) - \omega' (u \wedge t') | + | t - t' |
:= \| \omega_t - \omega'_{t'} \|_{\infty} + | t - t' |
\end{eqnarray*}
$ ( \Lambda_T, d_{\infty} ) $ is then a complete metric space.
Any map $F:[0, T] \times D( [0, T], \mathbb{R}^d)\to \mathbb{R}$ verifying the causality condition \eqref{eq:nonanticipative} can  be equivalently viewed as a functional on the space $\Lambda_T$ of stopped paths. 

We call 
{\it non-anticipative functional} on $ D( [0, T], \mathbb{R}^d ) $ a measurable map $F:(\Lambda_T,d_{\infty}) \to \mathbb{R}$ and denote 
 by $ \mathbb{C}^{0, 0} (\Lambda^d_T) $ the set of  such maps which are continuous for the $d_{\infty})$ metric. Some weaker notions of continuity for non-anticipative functionals turn out to be useful \cite{CF10A}:
\begin{definition} \label{def_functional}
A non-anticipative functional $ F $ is said to be:
\begin{itemize}
\item continuous at fixed times if for any $ t \in [0, T] $, $ F (t, \cdot): D( [0,T], \| \cdot \|_{\infty} )\to\mathbb{R}$ in $ [0, T] $, i.e. $ \forall \omega \in D( [0, T], \mathbb{R}^d ) $, $ \forall \epsilon > 0 $, $ \exists \eta > 0 $, $ \forall \omega' \in D( [0, T], \mathbb{R}^d ) $,
\[  \| \omega_t - \omega'_t \| _{\infty}< \eta \implies | F (t, \omega) - F (t, \omega') | < \epsilon \]
\item left-continuous if $ \forall (t, \omega) \in \Lambda_T $, $ \forall \epsilon > 0 $, $ \exists \eta > 0 $ such that $ \forall (t', \omega') \in \Lambda_T $,
\[ ( t' < t \quad \mathrm{and} \quad d_{\infty} \left( (t, \omega), (t', \omega') \right) < \eta ) \implies | F (t, \omega) - F (t', \omega') | < \epsilon \]
We denote by $ \mathbb{C}^{0, 0}_l (\Lambda^d_T) $ the set of left-continuous functionals. Similarly, we can define the set $ \mathbb{C}^{0, 0}_r (\Lambda^d_T) $ of right-continuous functionals.
\end{itemize}
\end{definition}
We also introduce a notion of local boundedness for functionals:
\begin{definition} \label{def_bound}
A non-anticipative functional $ F $ is said to be boundedness-preserving if for every compact subset $ K $ of $ \mathbb{R}^d $, $ \forall t_0 \in [0, T] $, $ \exists C (K, t_0) > 0 $ such that:
\[ \forall t \in [0, t_0],\quad \forall (t, \omega) \in \Lambda_T, \quad\omega ([0, t]) \subset K \implies | F (t, \omega) | < C (K, t_0). \]
We denote by $ \mathbb{B} (\Lambda^d_T) $ the set of boundedness-preserving functionals.
\end{definition}
We now recall some notions of differentiability for functionals following \cite{CF10B,cont2012,D09}. For $ e \in \mathbb{R}^d $ and $ \omega \in D ( [0,T], \mathbb{R}^d ) $, define the vertical perturbation $ \omega_t^e $ of $(t,\omega)$ as the  path obtained by shifting the path by $e$ at $t$:
$$\omega_t^e = \omega_t + e {\mathbf 1}_{[t, T]}.$$
\begin{definition} \label{def_diff}
A non-anticipative functional $ F $ is said to be:
\begin{itemize}
\item horizontally differentiable at $ (t, \omega) \in \Lambda_T $ if
\ba \mathcal{D} F (t, \omega) = \lim\limits_{h \to 0^+} \frac{F (t+h, \omega_t) - F (t, \omega_t)}{h} \label{eq.horizontal}\ea
exists. If $\mathcal{D} F (t, \omega)$ exists for all $ (t, \omega) \in \Lambda_T $, then \eqref{eq.horizontal} defines a non-anticipative functional $ \mathcal{D} F $, called the horizontal derivative of $ F $.
\item vertically differentiable at $ (t, \omega) \in \Lambda_T $ if the map:
\[
\begin{array}{ccc}
g_{(t,\omega)}:\mathbb{R}^d & \longrightarrow & \mathbb{R} \\
e & \mapsto & F (t, \omega_t + e {\mathbf 1}_{[t, T]})
\end{array}
 \]
is differentiable at  $ 0 $. Its gradient at $ 0 $ is called the Dupire derivative (or vertical derivative) of $ F $ at $ (t, \omega) $:
\ba \nabla_\omega F(t,\omega)= \nabla  g_{(t,\omega)} (0)\in \mathbb{R}^d\label{eq.vertical}\ea
i.e. $\nabla_{\omega} F (t, \omega) = ( \partial_i F (t, \omega), i = 1, \cdots, d) $
with
\[ \partial_i F (t, \omega) = \lim_{h \to 0} \frac{F (t, \omega_t + h e_i {\mathbf 1}_{[t, T]}) - F (t, \omega_t)}{h} \]
where $ ( e_i, i=1, \cdots, d ) $ is the canonical basis of $ \mathbb{R}^d $. If $ F $ is vertically differentiable at all $ (t, \omega) \in \Lambda^d_T $,   $ \nabla_{\omega} F: \Lambda_T\to \mathbb{R}^d$ defines a non-anticipative functional called the vertical derivative of $ F $.
\end{itemize}
\end{definition}
We may repeat the same operation on $ \nabla_{\omega} F $ and define similarly $ \nabla_{\omega}^2 F $, $ \nabla_{\omega}^3 F $, etc. 
\begin{remark}
Note that $\mathcal{D}F(t,\omega)$ is {\bf not} the partial derivative in $t$:
$$\mathcal{D}F(t,\omega)\neq \partial_tF(t,\omega)= \lim_{h\to 0}\frac {F(t+h, \omega)-F(t, \omega)} h.$$
$ \partial_tF(t,\omega)$ corresponds to a 'Lagrangian' derivative which follows the increment of $F$ along the path from $t$ to $t+h$, whereas in \eqref{eq.horizontal} the path is stopped at $t$ and the increment is taken along the stopped path.
\end{remark}
We consider the following classes of smooth functionals:
\begin{definition} \label{def:C12}
We define
$ \mathbb{C}^{1, k}_b (\Lambda^d_T) $ as the set of non-anticipative functionals
 $ F: (\Lambda^d_T,d_\infty) \to \mathbb{R}$ which are
\begin{itemize}
\item horizontally differentiable with $ \mathcal{D} F $ continuous at fixed times;
\item $ k $ times vertically differentiable with $ \nabla_{\omega}^j F \in \mathbb{C}^{0, 0}_l (\Lambda^d_T) $ for $ j = 0, \cdots, k $;
\item  $ \mathcal{D} F, \nabla_{\omega} F, \cdots, \nabla_{\omega}^k F \in \mathbb{B} (\Lambda^d_T) $.
\end{itemize}
\end{definition}


Consider now a sequence 
$\pi^n=(0=t^{n}_0<t^n_1<\cdots <t^{n}_{k(n)}=T)$ of partitions of $[0,T]$. 
$|\pi^n|=\sup\{ |t^n_{i+1}-t^n_i|, i=1\cdots k(n)\}\to 0$ will denote the mesh size of the partition.
A c\`adl\`ag path $x \in D([0,T],\mathbb{R})$ is said to have finite
quadratic variation along the sequence of partitions
$(\pi^n)_{n\geq 1}$ if for any $t\in [0,T]$ the limit
\ba [x]_\pi(t):=\mathop{\lim}_{n\to\infty}\sum_{t^n_{i+1}\leq t} (x(t^n_{i+1})-x(t^n_i))^2 <\infty\label{eq:qv}\ea
exists and the increasing function $[x]$ has  Lebesgue decomposition
$$[x]_\pi(t)=[x]_\pi^c(t) + \sum_{0 < s \leq t} |\Delta x(s)|^2$$
where $[x]_\pi^c$ is a continuous, increasing function. We denote the set of such paths $Q_\pi([0,T],\mathbb{R})$.

A d-dimensional path $x=(x^1,...,x^d) \in D([0,T],\mathbb{R}^d)$ is said to have finite
quadratic variation along $\pi$ if $x^i\in  Q_\pi([0,T],\mathbb{R})$ and $x^i+x^j\in  Q_\pi([0,T],\mathbb{R})$ for all $i,j=1..d$.
Then for any  $i,j=1,\cdots, d$  and $t\in [0,T]$, we have 
$$\sum_{t_k^n\in \pi^n, t^n_k\leq t}(x^i(t^n_{k+1})-x^i(t^n_k)).(x^j(t^n_{k+1})-x^j(t^n_k)) \mathop{\to}^{n\to\infty} [x]^{ij}_\pi(t)=\frac{[x^i+x^j]_\pi(t)-[x^i]_\pi(t)-[x^j]_\pi(t)}{2}.$$
The matrix-valued function $[x]:[0,T]\to S^+_d$  whose elements are given by $$[x]_\pi^{ij}(t)=\frac{[x^i+x^j]_\pi(t)-[x^i]_\pi(t)-[x^j]_\pi(t) }{2}$$
is called the {\it quadratic covariation} of the path $x$ along the sequence of partitions $\pi$. 
For further discussion of this concept, we refer to \cite{cont2012}.

Consider now a path $\omega \in Q_\pi([0,T],\mathbb{R}^d)$ with finite quadratic variation along $\pi$. Since $\omega$ has at most a countable set of jump times, we may  assume that the partition 'exhausts' the jump times in the sense that \ba \sup_{t \in
[0,T]-\pi^n} |\omega(t)-\omega(t-)|  \mathop{\to}^{n\to\infty} 0.\label{eq:exhaustion}\ea
Then the piecewise-constant approximation
\ba \omega^n(t)=\sum_{i=0}^{k(n)-1} \omega(t_{i+1}-)1_{[t_i,t_{i+1}[ }(t)+\omega(T)1_{\{T\}}(t) \label{piecewiseconstant.eq}\ea
 converges uniformly to $\omega$:
$\mathop{\sup}_{t\in [0,T]}\|\omega^n(t)-\omega(t)\|\mathop{\to}_{n\to\infty} 0.$
Note that with the notation \eqref{piecewiseconstant.eq}, $\omega^n(t^n_i-)=\omega(t^n_i-)$ but
$\omega^n(t^n_i)=\omega(t^n_{i+1}-).$ If we define
$$ \omega_{t^n_i}^{n, \Delta\omega(t_i^n)}= \omega^n +\Delta\omega(t^n_i) 1_{[t_i^n,T]},\quad{\rm then}\qquad
\omega^{n, \Delta\omega(t_i^n)}_{t^n_i-}(t^n_i)=\omega(t^n_{i}).$$

Approximating the variations of the functional by vertical and horizontal increments along the partition $\pi^n$, we obtain
the following pathwise change of variable formula for $\mathbb{C}^{1,2}(\Lambda^d_T)$ functionals, derived in \cite{CF10B} under more general assumptions:
\begin{theorem}[Pathwise change of variable formula for $\mathbb{C}^{1,2}$ functionals \cite{CF10B}] \label{thm.pathwise} 
Let $\omega\in Q_\pi([0,T],\mathbb{R}^d)$ verifying \eqref{eq:exhaustion}.
Then for any $F \in \mathbb{C}^{1,2}_{\rm b}(\Lambda^d_T)$ the limit
\ba \int_{0}^T \nabla_{\omega} F(t,\omega_{t-} ) d^\pi\omega
:= \mathop{\lim}_{n\to\infty} \sum_{i=0}^{k(n)-1} \nabla_{\omega}
F\left(t^n_{i},\omega^{n,\Delta \omega(t^n_i)
}_{t^n_i-}\right).
(\omega(t^n_{i+1})-\omega(t^n_i))\qquad\label{eq.Follmerintegral} \ea
exists and
\ba  F(T,\omega_T)-F(0,\omega_0) =\int_{0}^T
\mathcal{ D}F(t,\omega_{t} ) dt
+ \int_{0}^T \frac{1}{2} {\rm
tr}\left({}^t\nabla^2_{\omega} F(t,\omega_{t-})d[\omega]_\pi^c(t)\right) \nonumber \\
+ \int_{0}^T \nabla_{\omega} F(t,\omega_{t-} ) d^\pi\omega
+ \sum_{t \in ]0,T]} \left( F(t,\omega_t)-F(t,\omega_{t-}
)-\nabla_{\omega} F(t,\omega_{t-} ).\Delta \omega(t)\right).\nonumber\ea
\end{theorem}
A consequence of this theorem is the ability to define the pathwise integral  $\int_{0}^. \nabla_{\omega} F(t,\omega_{t-} ) d^\pi\omega $
 as a limit of "Riemann sums" computed along the sequence of partitions $\pi$. 
For a continuous path $\omega\in C^0([0,T],\mathbb{R}^d)$ this simplifies to:
\ba \int_{0}^T \nabla_{\omega} F(t,\omega ) d^\pi\omega
:= \mathop{\lim}_{n\to\infty} \sum_{i=0}^{k(n)-1} \nabla_{\omega}
F\left(t^n_{i},\omega^{n}_{t^n_i}\right).
(\omega(t^n_{i+1})-\omega(t^n_i))\label{eq.pathwiseintegral} \ea
This integral, first constructed in \cite{CF10B}, extends H. F\"ollmer's construction \cite{follmer1981} for integrands of the form $\nabla f,\ f\in C^2(\mathbb{R}^d)$ to (path-dependent) integrands of the form 
$$\nabla_\omega F,\qquad F\in \mathbb{C}^{1,2}_b(\Lambda^d_T).$$
This construction applies to typical paths of Brownian motion and semimartingales, which satisfy the quadratic variaton property almost-surely, although they have infinite $p$-variation for $p=2$ so \eqref{eq.pathwiseintegral} cannot be constructed as a Young integral \cite{young1936}.

The goal of this paper is to explore some properties of this integral.

\section{Isometry property of the pathwise integral}\label{sec:isometry}
\subsection{Pathwise isometry formula}
 
Let $\omega\in Q_\pi([0,T],\mathbb{R}^d)$ be a given path with finite quadratic variation along a nested sequence of partitions $\pi=(\pi^n)_{n\geq 1}$.
Our goal is to provide conditions under which the the following identity holds:
		\begin{equation}\label{isometry1forms}
		\left[F(\cdot,\omega_{\cdot})\right]_\pi(t)=\int\limits_0^t
		\langle \nabla_{\omega}  F(s,{\omega}_{s-})\ ^t\nabla_{\omega}  F (s,{\omega}_{s-}), d[{\omega}]_\pi(s)\rangle,
		\end{equation}
		where we use the notation $\langle A , B\rangle :=tr(A B)$, for square matrices $A, B.$
		
This relation was first noted in \cite{cont2012} for functionals of the form $F(\omega)=\int_0^t\phi.d\omega$ with $\phi$ piecewise-constant.
		To extend this property to a more general setting we assume a Lipschitz-continuity condition on the functional $F$ and a H\"older-type  regularity condition for the path $\omega$.
\begin{hyp}[Uniform Lipschitz continuity]\label{ass.Fholder}
	$$\exists K>0, \quad \forall \omega,\, \omega'\in D([0,T], \Rr^d), \forall t\in [0,T],\quad
		|F(t,\omega)-F(t,\omega')|\leq K \|\omega_t-\omega'_t\|_{\infty}.$$
	
\end{hyp}
We denote by $Lip(\Lambda^d_T, \|\cdot\|_{\infty})$ 	 the space of all non-anticipative functionals satisfying the above property.

Consider a nested sequence of partitions $\pi^n=\{t^n_i,i=0, \ldots, m(n) \}$ of $[0,T]$. 
We define the oscillation  osc$( f, \pi^n)$ of a function $f:[0,T]\to\mathbb{R}^d$ between points of the partition $\pi^n$ as
	\[
	{\rm osc}(f, \pi^n)=
\mathop{\max}_{j=0..m(n)-1} \mathop{\sup}_{t \in ( t^n_j, t^n_{j+1}] } | f(t) - f( t^n_j) |
	\]
We make the following assumptions on the sequence of partitions:
\begin{hyp}[Vanishing oscillation along $\pi$] \label{ass.oscill}
$${\rm osc}(\omega, \pi^n)\mathop{\to}^{n\to+\infty}  0.$$
\end{hyp}
\begin{hyp}\label{ass.Foscill}
 $$\max_{i\in \{0, \ldots, \ m(n)-1\}}|F(t^n_{i+1},\omega)-F(t^n_i, \omega)|\mathop{\to}^{n\to+\infty} 0.$$
\end{hyp}
\begin{remark}[Partitions with vanishing mesh]
If $|\pi^n|\to 0$ and $F:(\Lambda^d_T,d_\infty)\to \mathbb{R}$ is Lipschitz-continuous then Assumptions \ref{ass.oscill} and \ref{ass.Foscill} are satisfied for any $\omega\in  C^0([0,T],\mathbb{R}^d)$.  

To see this, set $\epsilon >0$. By uniform continuity of $\omega$  on $[0,T]$  there exists $\delta(\epsilon) >0$ such that $\|\omega(t)-\omega(s)\|\leq \epsilon$ whenever $|t-s|\leq \delta$.
Since $|\pi^n|\to 0$ there exists $N_1\leq 1$ such that  $|\pi^n|\leq \epsilon \wedge \delta(\epsilon)$   for $n\leq N_1$. Then for $n\geq N_1$ we have
$ {\rm osc}(\omega, \pi^n)\leq \epsilon $ so Assumption \ref{ass.oscill}  holds. Furthermore, for $n\leq N_1$,
$$d_\infty( (t^n_{i+1},\omega),  (t^n_i, \omega))=|t^n_{i+1}-t^n_{i}| + \|\omega_{t^n_{i+1}}-\omega_{t^n_{i}} \|_\infty$$
$$\leq |\pi^n|+ \mathop{\sup}_{t\in [t^n_{i},t^n_{i+1} ]}|\omega(t)-\omega(t^n_{i})| \leq 2\epsilon $$
so $d_\infty( (t^n_{i+1},\omega),  (t^n_i, \omega))\to 0$ as $n\to\infty$.  Lipschitz-continuity of $F$ then entails $\max_{i\in \{0, \ldots, \ m(n)-1\}}|F(t^n_{i+1},\omega)-F(t^n_i, \omega)|\to^{n\to+\infty} 0.$
\end{remark}
In section \ref{sec:LebesguePartition} we will also consider sequences of partitions whose mesh may not converge to zero. For such partitions
Assumption \ref{ass.Foscill} is not redundant.

For $0< \nu < 1,$ denote by $C^{\nu}([0,T])$ the space of $\nu-$H\"older continuous functions:
$$ C^{\nu}([0,T],\mathbb{R}^d) =\{ f\in C^0([0,T]),\quad \sup_{(t,s)\in[0,T]^2, t\neq s }\frac{\|f(t)-f(s)\|}{|t-s|^{\nu}}<+\infty \},$$ $${\rm and}\qquad   C^{\nu-}([0,T],\mathbb{R}^d) =\mathop{\cap}_{0\leq \alpha< \nu}C^{\alpha}([0,T],\mathbb{R}^d) $$ the space of functions which are $\alpha-$H\"older for every $\alpha<\nu$.
For $\omega\in C^\nu([0,T],\mathbb{R}^d),$ the following piecewise constant approximation
\begin{equation}\label{stepbypath}
	\omega^n:=\sum_{i=0}^{m(n)-1}\omega(t^n_{i+1}-)\mathbf{1}_{[t^n_i, t^n_{i+1})}+\omega(T)\mathbf{1}_{\{T\}}
\end{equation}
	satisfies: $\|\omega-\omega^n\|_\infty \leq |\pi^n|^{\nu}$.

Our main result  is the following:
\begin{theorem}[Pathwise Isometry formula] 	Let $\pi=(\pi^n)_{n\geq 1}$ be a sequence of partitions of $[0,T]$, and $\omega\in Q_\pi([0,T],\mathbb{R}^d)\cap C^{\nu}([0,T],\mathbb{R}^d)$ for $\nu > \frac{\sqrt{3}-1}{2}$, satisfying Assumption \ref{ass.oscill}.

Let $F \in \mathbb{C}^{1,2}_b(\Lambda^d_T)\cap Lip(\Lambda^d_T, \|\cdot\|_{\infty})$  with  $\nabla_{\omega} F\in \mathbb{C}^{1,1}_b(\Lambda^d_T)$  such that
Assumption  \ref{ass.Foscill} holds. Then:
	\begin{equation}\label{eq.isometry1forms}
	\left[F(\cdot,\omega_{\cdot})\right]_\pi (T)= \left[ \int_0^.
	\nabla_{\omega}  F (s,{\omega}_{s}) d^\pi\omega\right]_\pi(T)=\int\limits_0^T
	\langle\nabla_{\omega}  F (s,{\omega}_{s})\ ^t\nabla_{\omega}  F (s,{\omega}_{s}) ,  d[{\omega}]_\pi \rangle.
	\end{equation}	\label{theorem.isometry}
\end{theorem}

\begin{remark} We note that for typical paths of a Brownian diffusion  or continuous semimartingale with non-degenerate local martingale component, the assumptions of Theorem \ref{theorem.isometry} hold almost surely as soon as Assumption \ref{ass.Fholder} is satisfied.
	More generally the result holds for $\omega\in C^{\frac{1}{2}-}([0,T],\mathbb{R}^d)$, which corresponds to the H\"older regularity of Brownian paths \cite{revuzyor}.
\end{remark}

We start with an important lemma  which is essential in proving both the  pathwise isometry formula and the uniqueness of pathwise integral. 
This result also establishes  a connection between the  functional derivatives introduced in Section \ref{sec:FIC} and controlled rough paths \cite{gubinelli2004}.

\begin{lemma}\label{theorem.controll} Let $\omega\in C^{\nu}([0,T],\mathbb{R}^d)$ for some ${\nu}\in(1/3, 1/2]$ and $F\in \mathbb{C}^{1,2}_b(\Lambda^d_T,\mathbb{R}^n)$  with $ \nabla_{\omega} F\in \mathbb{C}^{1,1}_b(\Lambda^d_T, \mathbb{R}^{n\times d})$  and $ F\in Lip(\Lambda^d_T, \|\cdot\|_{\infty})$. 
	Define
	\begin{equation}\label{eq.controlledpath}
		R^{F}_{s,t}(\omega):=F(s,\omega_s)-F(t,\omega_t)-\nabla_{\omega}F(t,\omega_t)(\omega(s)-\omega(t)).
	\end{equation}	
	Then there exists a constant $C_{F,T,\|\omega\|_{\nu}}$ increasing in $T$ and $\|\omega\|_{\nu}$, such that
	$$|R^{F}_{s,t}(\omega) |\leq C_{F,T,\|\omega\|_{\nu}}  |s-t|^{2\beta_{\nu}},$$
	with $2 \beta_{\nu}=\nu(1+\nu)$. 
	
\end{lemma}

\begin{proof}
		We will prove only the case when the values of $F$ are scalar, i.e. $n=1$; the extension to the general case is straightforward.
		
		We start by recalling the following result from  \cite[Proposition 5.26]{cont2012}:
		\begin{lemma}
			Let $G\in \mathbb{C}^{1,1}_b(\Lambda^d_T)$ and $\lambda:[t,s]\to\mathbb{R}$  a continuous function with finite variation. Then
			\begin{equation*}
			G(s,\lambda_s)-G(t,\lambda_t)=\int_t^s\mathcal{D} G(u,\lambda_u)du+\int_t^s\nabla_{\omega} G(u,\lambda_u)d\lambda(u).\tag{*}
			\end{equation*}
		Here  the integrals are defined in the Riemann-Stieltjes sense.
\end{lemma}
		
		Now, let us fix a  Lipschitz continuous  path $\lambda\colon [0,T]\to\Rr^d$, using the above lemma repeatedly for $G=F,\, \nabla_{\omega} F, \cdots$ we will obtain an integral formula for $R^{F}_{t,s}(\lambda)$ in terms of the derivatives $F$. For the sake of convenience we denote by $\partial_i F$ and $\lambda^i$ respectively the $i$-th coordinates of $\nabla_{\omega} F$ and $\lambda.$ We also use Einstein's convention of summation in repeated indexes and the following notation
		\[
		\delta \lambda_{t,s}:=\lambda(s)-\lambda(t).
		\]
		
		Using (*) for $G=F$, we have
		
		\begin{equation}\label{eq.Rfy1}
		R^{F}_{t,s}(\lambda)=\int_t^s\mathcal{D}F(u,\lambda_u)du+\int_t^s\left(\partial_i F(u,\lambda_u)-\partial_i F(t, \lambda_t)\right)d\lambda^i(u)
		\end{equation}
		For the second term  on the right-hand side of the above identity we use the (*) with $G=\partial_{i} F$  and then Fubini's theorem to get
		\begin{equation}\label{eq.Rfy2}
		\begin{split}
		\int_t^s\left(\partial_i F(u,\lambda_u)-\partial_i F(t, \lambda_t)\right)d\lambda^i(u)=\int_t^s\int_t^u \mathcal{D}\partial_i F(r,\lambda_r)drd\lambda^i(u)\\
		+ \int_t^s\int_t^u \partial^2_{ij} F(r,\lambda_r)d\lambda^j(r)d\lambda^i(u)\\
		=\int_t^s\mathcal{D}\partial_i F(r,\lambda_r)(\lambda^i(s)-\lambda^i(r))dr +\int_t^s\partial^2_{ij} F(r,\lambda_r)\dot{\lambda}^j(r)(\lambda^i(s)-\lambda^i(r))dr.
		\end{split}
		\end{equation}
		
		Thus
		\begin{equation}\label{eq.Rfyformula}		
		R^{F}_{t,s}(\lambda)=\int_t^s ( \mathcal{D}F(u,\lambda_u)+\mathcal{D}\partial_i F(u,\lambda_u)\delta\lambda^i_{u,s} )du
		+\int_t^s\partial^2_{ij} F(r,\lambda_r) \delta\lambda^i_{r,s}d\lambda^j(r)		
		\end{equation}
		
To use the above formula for estimate the error term $R^{F}$  for H\"older continuous  paths we will  use a piecewise linear approximation.	For $\omega\in C^{\nu}$ let $\omega^N$ be the piecewise linear approximation of $\omega$  on $[t,s]$ defined by
		\begin{itemize}
			\item $\omega^N(r)=\omega(r),\quad \forall r\in[0,t]$,\\
			
			\item $\omega^N(\tau^N_k)=\omega(\tau^N_k),\,k=\{0,...,N\}$, where  $\tau^N_k=t+ k\frac{s-t}N$,\\
			
			\item $\omega^N$ is linear on each interval $[\tau^N_k,\tau^N_{k+1}]$.
		\end{itemize}
		Then, we have
		
		\begin{align}\label{eq.omgaNbounds}
		\|\omega^N-\omega\|_{\infty}\leq C \|\omega\|_{\nu} \frac {|s-t|^{\nu}}{N^{\nu}},\quad \|\omega^N\|_{\nu}\leq C\|\omega\|_{\nu},\nonumber
		\\
		|\dot{\omega}^N|=\sum_i \frac{|\delta\omega_{\tau_i,\tau_{i+1}}|}{|\tau_{i+1}-\tau_i|} 1_{(\tau_i,\tau_{i+1})}\leq N^{1-\nu} \|\omega\|_{\nu}|s-t|^{\nu-1},
		\\
		|\delta\omega_{a,b}|\leq \|\omega\|_{\nu}|b-a|^{\nu}.\nonumber
		\end{align}

		Using the local boundedness property of $\mathcal{D}F, \mathcal{D}\nabla_{\omega} F$ and $\nabla_{\omega} ^2 F$, and the representation \eqref{eq.Rfyformula} for $\lambda=\omega^N$,
		we obtain
		\begin{equation}\label{eq.Rfylongestimate}
		\begin{split}
		|R^{F}_{t,s}(\omega^N)|\leq C_F|s-t|+
		C_F\|\omega^N\|_{\nu}|s-t|^{1+\nu} + C_F\|\omega^N\|^2_{\nu} N^{1-\nu} |s-t|^{2\nu}.
		\end{split}
		\end{equation}

		On the other hand since $\omega^N_t=\omega_t$ and $\omega^N(s)=\omega(s)$
		\begin{equation}\label{eq.Rfyapprox}
		|R^{F}_{t,s}(\omega^N)-R^{F}_{t,s}(\omega)|=|F(s,\omega^N)-F(s,\omega)|\leq C_F\|\omega^N-\omega\|_{\infty}\leq C_F\|\omega\|_{\nu}N^{-\nu}|s-t|^{\nu}.
		\end{equation}
		
		from above estimates \eqref{eq.Rfylongestimate}, \eqref{eq.Rfyapprox} and triangle inequality
		\[
		|R^{F}_{t,s}(\omega)|\leq C_F|s-t|+
		C_F\|\omega\|_{\nu}|s-t|^{1+\nu}+C_F\|\omega\|^2_{\nu} N^{1-\nu} |s-t|^{2\nu}+N^{-\nu}|s-t|^{\nu}.
		\]
		To optimize the above bound, we choose $N$ so that
		$
		\|\omega\|^2_{\nu}  N^{1-\nu}|s-t|^{2\nu}\approx N^{-\nu}|s-t|^{\nu}
		$
		i.e.
		$
		N\approx\|\omega\|_{\nu} ^{-2} |s-t|^{-\nu}.
		$
		Hence the result
		\[
		|R^{F}_{t,s}(\omega)|\leq C_{F,T}(1+	\|\omega\|_{\nu})|s-t| +C_{F,\nu}
		\|\omega\|_{\nu}^{2\nu}|s-t|^{\nu+\nu^2}.
		\]

\end{proof}

We are now ready to prove the isometry property:
\begin{proof}[\textbf{Proof of  Theorem \ref{theorem.isometry}}] To prove the isometry formula we note that 
	
	\[ 
	\begin{split}
	\left|	\left(F(t^n_{i+1},\omega_{t^n_{i+1}})-F(t^n_i, \omega_{t^n_i})\right)^2-\langle \nabla_{\omega}F(t,\omega_{t^n_i})\, ^t \nabla_{\omega}F(t,\omega_{t^n_i}),  \delta\omega_{t^n_i, t^n_{i+1}}\, ^t\delta\omega_{t^n_i, t^n_{i+1}}\rangle\right|\\
	\leq  |R^{F}_{t^n_i, t^n_{i+1}}(\omega)|^2 +C_F |R^{F}_{t^n_i,t^n_{i+1}}(\omega)| |\delta\omega_{t^n_i, t^n_{i+1}}|
	\end{split}
	\] 
	From the Assumption \ref{ass.oscill} it follows that $M_n:=\max_i\left| R^{F}_{t^n_i, t^n_{i+1}}(\omega)\right|\to_{n\to\infty} 0$ and from part a) of Theorem \eqref{theorem.controll} we have $\left| R^{F}_{t^n_i, t^n_{i+1}}(\omega)\right|\leq C |t^n_{i+1}-t^n_i|^{\nu^2+\nu}$. Since  $\nu>(\sqrt{3}-1)/2$ is equivalent to $\nu^2+\nu>\frac 1 2$, we get
	\[
	\sum_i\left| R^{F}_{t^n_i, t^n_{i+1}}(\omega)\right|^2\leq C M_n^{2-\frac{1}{\nu^2+\nu}}\sum_i|t^n_{i+1}-t^n_i|\leq C TM_n^{2-\frac{1}{\nu^2+\nu}}\to 0.
	\]
	Consequently, since $\sum_i|\delta\omega_{t^n_i, t^n_{i+1}}|^2\to tr([\omega]^{\pi})$, we have by Cauchy -Schwarz inequality
	\[
	\sum_i|R^{F}_{t^n_i,t^n_{i+1}}(\omega)| |\delta\omega_{t^n_i, t^n_{i+1}}|\leq \sqrt{	\sum_i\left| R^{F}_{t^n_i, t^n_{i+1}}(\omega)\right|^2} \sqrt{	\sum_i|\delta\omega_{t^n_i, t^n_{i+1}}|^2}\to 0.
	\]
	We conclude
	\[ 
	\begin{split}
	\left|	\sum_i\left(F(t^n_{i+1},\omega_{t^n_{i+1}})-F(t^n_i, \omega_{t^n_i})\right)^2-\sum_i \langle \nabla_{\omega}F(t,\omega_{t^n_i})\, ^t \nabla_{\omega}F(t,\omega_{t^n_i}),  \delta\omega_{t^n_i, t^n_{i+1}}\, ^t\delta\omega_{t^n_i, t^n_{i+1}}\rangle\right|\\
	\leq  \sum_i|R^{F}_{t^n_i, t^n_{i+1}}(\omega)|^2 +C_F \sum_i|R^{F}_{t^n_i,t^n_{i+1}}(\omega) | |\delta\omega_{t^n_i, t^n_{i+1}}|\to 0
	\end{split}
	\] 
	which together with 
	\[
	\sum_i \langle \nabla_{\omega}F(t,\omega_{t^n_i})\, ^t \nabla_{\omega}F(t,\omega_{t^n_i}),  \delta\omega_{t^n_i, t^n_{i+1}}\, ^t\delta\omega_{t^n_i, t^n_{i+1}}\rangle\to\int\limits_0^T
	\langle \nabla_{\omega}  F (s,{\omega}_{s})\ ^t\nabla_{\omega}  F (s,{\omega}_{s}) ,  d[{\omega}]_\pi \rangle.
	\]
	 implies 
	\[
	 \left[ \int_0^.
	\nabla_{\omega}  F (s,{\omega}_{s}) d^\pi\omega\right]_\pi(T)=\int\limits_0^T
	\langle \nabla_{\omega}  F (s,{\omega}_{s})\ ^t\nabla_{\omega}  F (s,{\omega}_{s}),  d[{\omega}]_\pi \rangle.
	\]

To prove the  equality $\left[F(\cdot,\omega_{\cdot})\right]_\pi(T)= \left[ \int_0^.
\nabla_{\omega}  F (s,{\omega}_{s}) d^\pi\omega\right]_\pi(T)$	we note, using the change of variable formula (Theorem \eqref{thm.pathwise})   it is enough to prove that the paths $\int_{0}^{\cdot}
\mathcal{ D}F(t,\omega_{t} ) dt$, $\int_{0}^{\cdot} \left<\nabla^2_{\omega} F(t,\omega_{t-}),d[\omega]_\pi(t)\right>$ have zero quadratic variation along $\pi^n.$ For the second path it follows from Assumption \ref{ass.oscill} as it implies that $\max_{i} |[\omega](t^n_{i+1})-[\omega]({t^n_i})|\to 0.$ For the first path note that since it is obviously a path of bounded variation we just need to prove that $\max_i \left|\int_{t^n_i}^{t^n_{i+1}}
\mathcal{ D}F(t,\omega_{t} ) dt \right|\to 0$.  For that we write
\[\begin{split}
\int_{t^n_i}^{t^n_{i+1}}
\mathcal{ D}F(t,\omega_{t} ) dt =\int_{t^n_i}^{t^n_{i+1}}
[\mathcal{ D}F(t,\omega_{t} )- \mathcal{ D}F(t,\omega_{t^n_i} )]dt +\int_{t^n_i}^{t^n_{i+1}}
\mathcal{ D}F(t,\omega_{t^n_i} ) dt \\
=\int_{t^n_i}^{t^n_{i+1}}
[\mathcal{ D}F(t,\omega_{t} )- \mathcal{ D}F(t,\omega_{t^n_i} )]dt +F(t^n_{i+1}, \omega_{t^n_i} ) -F(t^n_i, \omega_{t^n_i} )\\
=\int_{t^n_i}^{t^n_{i+1}}
[\mathcal{ D}F(t,\omega_{t} )- \mathcal{ D}F(t,\omega_{t^n_i} )]dt +\left(F(t^n_{i+1}, \omega_{t^n_i} )- F(t^n_{i+1}, \omega_{t^n_{i+1}}) \right)\\
+  \left(F(t^n_{i+1}, \omega_{t^n_{i+1}})  -F(t^n_i, \omega_{t^n_i} )\right).
\end{split}
\]
By the continuity of $F, \mathcal{ D}F$ and Assumption \ref{ass.oscill} the first two summands in the right hand side of above equality uniformly converge to $0$ as $n\to \infty$, by Assumption \ref{ass.Foscill} the third summand also converges uniformly to $0$. 
\end{proof}


\begin{remark}[Relation with Ito isometry formula]\label{rem:Itoisometry}
	Let $\mathbb{P}$ be the Wiener measure on $C^0([0,T],\mathbb{R})$. Then the pathwise integral
 $(t,\omega)\to \int_0^t \nabla_\omega F(u,\omega) d^\pi\omega$ is defined $\mathbb{P}-$almost surely and defines a process which is a version of the Ito integral $\int_0^.\nabla_\omega F(t,W)dW $. Integrating \eqref{eq.isometry1forms} with respect to $\mathbb{P}$ yields the well-known Ito isometry formula \cite{ikedawatanabe}:
$$  E^\mathbb{P}\left( [\int_0^.\nabla_\omega F(t,W)dW ](T) \right)=E^\mathbb{P}\left( |\int_0^T\nabla_\omega F(t,W)dW|^2\right)=E^\mathbb{P}\left( \int_0^T|\nabla_\omega F(t,W)|^2 dt \right). $$
So, Theorem \ref{theorem.isometry} reveals that underlying the Ito isometry formula is the pathwise isometry \eqref{eq.isometry1forms} which does not  rely on the Wiener measure. 
\end{remark}

\subsection{Partitions of hitting times}\label{sec:LebesguePartition}

 An alternative approach  to the definition of the pathwise integral is to consider quadratic variation and Riemann sums computed along 'Lebesgue partitions' i.e. partitions of hitting times of uniformly spaced levels \cite{chacon1981,davis2014,karandikar1995,perkowski2015}:
\ba \tau^n_0(\omega)=0,\qquad  \tau^n_{k+1}(\omega)= \inf\{ t>\tau^n_k, |\omega(t)-\omega(\tau^n_k)|\geq 2^{-n}\}\wedge T. \label{eq.hittingtimes}\ea
Although this is different in spirit from the original approach of F\"ollmer, there are some arguments for using such 'intrinsic' sequences of partitions for a given path, one of them being that  quadratic variation computed along  \eqref{eq.hittingtimes} is invariant under a time change. Thus a natural question is whether the results above apply to such partitions.
We will now present a version of Theorem \ref{theorem.isometry} for the sequence of partitions defined by \eqref{eq.hittingtimes}. 

Define  $m(n)= \inf\{k \geq 1, \tau^n_k(\omega)=T  \}$ and denote  $\tau^n(\omega)=(\tau_k^n(\omega), k=0,\ldots, m(n))$.
It is easy to see that the vanishing oscillation conditon (Assumption  \ref{ass.oscill})  automatically holds for the sequence $(\tau^n(\omega),n\geq 1)$ since by construction 
$${\rm osc}(\omega, \tau^n(\omega))\leq 2^{-n}.$$
However, as the example of a constant path shows, 
 for a general path $\omega\in C^0([0,T],\mathbb{R}^d)$, the partitions $\tau^n(\omega)$ may fail to have a vanishing mesh size. 

We will show now that, if a path $\omega$ has  {\it strictly increasing } quadratic variation along its Lebesgue partition $\tau^n(\omega)$ then $|\tau^n|\to 0$.

\begin{lemma}\label{lem.lebpart} Assume $\omega\in C^0([0,T],\mathbb{R}^d)$ has finite quadratic variation along the Lebesgue partition $\tau^n(\omega)$ defined by  \eqref{eq.hittingtimes}:
\ba\label{lebvariation} \forall t\in (0,T],\qquad \mathop{\lim}_{n\to\infty}\sum_{\tau^n_k\leq t} |\omega( \tau^n_{k+1})-\omega(\tau^n_k)|^2 =[\omega]_{\tau(\omega)}(t)>0.\label{eq.qvLebesgue}\ea
If $t\to [\omega]_{\tau(\omega)}(t)$ is a strictly increasing function
then 
 $|\tau^n(\omega)|\to 0$.
\end{lemma}

\begin{proof} 
Let $ h>0, t\in [0,T-h]$. Denote by $m(n,t,t+h)$  the number of partition points of $\tau^n(\omega)$ in the interval $[t,t+h]$. Then 
$$  \sum_{t\leq \tau^n_k\leq t+h} |\omega( \tau^n_{k+1})-\omega(\tau^n_k)|^2 =4^{-n} (m(n,t,t+h)-1).$$
Since $t\to [\omega]_{\tau(\omega)}(t)$ is strictly increasing,
$$
\sum_{t\leq \tau^n_k\leq t+h} |\omega( \tau^n_{k+1})-\omega(\tau^n_k)|^2 =4^{-n} (m(n,t,t+h)-1) \mathop{\to}_{n\to\infty}[\omega]_{\tau(\omega)}(t+h)-[\omega]_{\tau(\omega)}(t) >0$$
so in particular $m(n,t,t+h) \sim 4^n \to \infty$. This holds for any interval $t\in [0,T-h]$ which implies that the  $|\tau^n|< h$ for $n$ large enough. Since this is true for any $h>0$, this implies $|\tau^n|\to 0$.

\end{proof}
The  condition that $[\omega]_{\tau(\omega)}$ is strictly increasing is an 'irregularity' condition on $\omega$: it means that the quadratic variation over any interval is non-zero.

For paths verifying this condition,
Theorem \ref{theorem.isometry} may be extended to Lebesgue partitions:
\begin{theorem}[Pathwise Isometry formula along Lebesgue partitions] Let $\omega\in Q_{\tau(\omega)}([0,T],\mathbb{R}^d)\cap C^{\nu}([0,T],\mathbb{R}^d)$   for some $\nu \in \left(\frac{\sqrt{3}-1}{2}, \frac 1 2\right)$  and  $F \in \mathbb{C}^{1,1}_b(\Lambda^d_T)\cap Lip(\Lambda^d_T, \|\cdot\|_{\infty})$ such that $\nabla_{\omega} F\in \mathbb{C}^{1,1}_b(\Lambda^d_T)$. If $\omega$ has strictly increasing quadratic variation $[\omega]_{\tau(\omega)}$ with respect to $\tau(\omega)$ then the isometry formula holds along the partition of hitting times defined by \eqref{eq.hittingtimes}:
	\[
	\left[F(\cdot,\omega_{\cdot})\right]_{\tau(\omega)}(t)= \left[ \int_0^.
	\nabla_{\omega}  F (s,{\omega}_{s-}) d^\pi\omega\right](t)=\int\limits_0^t
	\langle \nabla_{\omega}  F (s,{\omega}_{s})\ ^t\nabla_{\omega}  F (s,{\omega}_{s}),  d[{\omega}]_{\tau(\omega)} \rangle.
	\]
\end{theorem}

\begin{proof}
First note that the sequence of partitions	$(\tau^n(\omega),n\geq 1)$  automatically satisfies  Assumption \ref{ass.oscill} since ${\rm osc}(\omega,\tau_n(\omega))\leq 2^{-n}.$

Since $[\omega]_{\tau(\omega)}$ is strictly increasing, Lemma \ref{lem.lebpart} implies that $|\tau^n|\to 0$. Since $\omega\in C^{\nu}([0,T],\mathbb{R}^d)$ which together with the continuity of $F$ implies that
$$\mathop{\max}_{i\in \{0, \ldots, \ m(n)-1\}}|F(t^n_{i+1},\omega)-F(t^n_i, \omega)|  \mathop{\to}^{n\to+\infty} 0.$$ so the sequence of partitions $\tau^n(\omega)$ satisfies Assumption  \ref{ass.Foscill}.

The result then follows from the Theorem \ref{theorem.isometry}.

\end{proof}

\subsection{The pathwise integral as a continuous isometry}\label{sec:quotient}
Theorem \ref{theorem.isometry}  suggests that the existence of an  isometric mapping underlying the pathwise integral. We will now proceed to make this structure explicit.

A consequence of Theorem \ref{theorem.isometry} is that, if $\omega$ is an irregular path  with strictly increasing quadratic variation, any regular functional $F(.,\omega)$ has zero quadratic variation along  $\pi$ if and only if  $\nabla_\omega F$ vanishes along $\omega$: 
\begin{proposition}[Preservation of irregularity]\label{0quadratic}
	Let $\omega\in  Q_\pi([0,T],\mathbb{R}^d)\cap C^{\nu}([0,T], \mathbb{R}^d)$ for some $\nu \in \left( \frac{\sqrt{3}-1}{2}, \frac 1 2\right]$ such that $\frac{d[\omega]^\pi}{dt}:=a(t) >0$  is a right-continuous, positive definite function, and $F \in \mathbb{C}^{1,2}_b(\Lambda^d_T)\cap Lip(\Lambda^d_T, \|\cdot\|_{\infty})$ be a non-anticipative functional with $ \nabla_{\omega} F\in \mathbb{C}^{1,1}_b(\Lambda^d_T)$.

If Assumptions \ref{ass.oscill} and \ref{ass.Foscill} hold, then the path $t\mapsto F(t,\omega)$ has  zero quadratic variation along the partition $\pi$ if and only if $\nabla_\omega F(t, \omega)=0,\, \forall t\in[0,T]$.
\end{proposition}
\begin{proof} Indeed, from Theorem \ref{theorem.isometry}
	\[
	\left[F(\cdot,\omega_{\cdot})\right](T)=\int\limits_0^T
	{}^t\nabla_{\omega}  F (s,{\omega}) a(s)\nabla_{\omega}  F (s,{\omega})\, ds.
	\]

	Since $a(t)$ is positive definite the integrand on the right hand side is non-negative and strictly positive unless $\nabla_{\omega}  F (s,{\omega}_{s})=0$. So by right-continuity of the integrand the integral is zero if and only if  $\nabla_\omega F(\cdot, \omega_{\cdot})\equiv 0$.
\end{proof}
Let $F\in C^{1,2}_b(\Lambda_T)$ be a functional satisfying assumptions of Theorem \ref{theorem.isometry}. We denote $\mathcal{R}(\Lambda_T)$ the set of such  {\it regular } functionals.

Proposition \ref{0quadratic} may be understood in the following way: if a path is  'irregular' in the sense of having  strictly increasing quadratic variation, this property is locally preserved by any regular functional transformation $F$  as long as $\nabla_{\omega}  F$ does not vanish.

We note that for $\omega\in C^{\frac{1}{2}-}([0,T],\mathbb{R}^d)$, we have $F(\cdot, \omega)\in C^{\frac{1}{2}-}([0,T],\mathbb{R})$. 

Let $a\colon [0,T]\to S_{+}^d$ be a continuous function taking values in positive-definite symmetric matrices. 
\begin{definition}[Harmonic functionals]  $F\in\mathcal{R}(\Lambda_T)$ is called $a-$harmonic   if
	\[
\forall (t,\omega)\in \Lambda_T,\qquad	  \mathcal{D} F(t, \omega_t) +\frac 1 2 \langle \nabla_{\omega}^2 F(t, \omega_t) , a(t) \rangle=0.
	\]
	We denote $\mathcal{H}_a(\Lambda_T)$  the space of $a(\cdot)$-harmonic functionals.
\end{definition}
Note that for any $F\in \mathcal{R}(\Lambda_T)$ the functional  defined by $$G(t,\omega)=F(t,\omega)-\int_0^t  \mathcal{D} F(s, \omega) ds-\frac 1 2 \int_0^t \langle \nabla_{\omega}^2 F(s, \omega) , a(s) \rangle ds$$ is $a(\cdot)$-harmonic.

The class of harmonic functionals plays a key role in probabilistic interpretation of the Functional Ito calculus \cite{fournie2010,ContFournie2013,cont2012}. We will now see that this class also plays a role in the extension of the pathwise integral.

Let $\bar\omega\in Q_\pi([0,T],\mathbb{R}^d)\cap C^{\nu}([0,T],\mathbb{R}^d)$ such that $d[\bar\omega]_\pi/dt = a$.
The functional change of variable formula (Theorem \ref{thm.pathwise} ) then implies that
$$\forall F\in \mathcal{H}_a(\Lambda_T),\quad \forall t\in [0,T], \qquad F(t,\bar\omega)=F(0,\bar\omega)+\int_0^t \nabla_\omega F(u,\bar\omega)d^\pi\bar\omega.$$
By Theorem \ref{theorem.isometry}, we have 
$$ [F(.,\bar\omega)]_\pi(t) = \int_0^t {}^t\nabla_\omega F(u,\bar\omega) a(u) \nabla_\omega F(u,\bar\omega) du= \|\nabla_\omega F(.,\bar\omega)\|^2_{L^2([0,T],a)} < \infty.$$
Consider the vector spaces
\ba
\begin{split}
\mathcal{H}_a(\bar\omega)&:=\left\{\,F(\cdot, \bar\omega_{\cdot})\,\big|\, F\in \mathcal{H}_a(\Lambda_T) \right\}\subset Q_\pi([0,T],\mathbb{R}),\\
\mathbb{V}_a(\bar\omega)&:=\left\{\,\nabla_{\omega}F(\cdot, \bar\omega_{\cdot})\,\big|\, F\in \mathcal{H}_a(\Lambda_T) \right\} \subset L^2([0,T],a).\nonumber
\end{split}
\ea
which are the images of $\overline{\omega}$ under all harmonic functionals and their vertical derivatives.
Proposition \ref{0quadratic} then implies that the map
\ba 
\omega\in \mathcal{H}_a(\bar{\omega})&\to   & \|\omega\|_\pi=\sqrt{[\omega]_\pi(T)}\nonumber\ea
defines a norm on $\mathcal{H}_a(\bar{\omega})$.
The pathwise integral can then be lifted to a continuous map
\ba
I_{\bar\omega}\colon \left(\mathbb{V}_a(\bar\omega), \|\cdot\|_{L^2([0,T], a)}\right)& \to &  \left(\mathcal{H}_a(\bar{\omega}), \|\cdot\|_{\pi}\right)\nonumber\\
\phi= \nabla_\omega F(.,\bar\omega ) & \to & I_{\bar\omega}(\phi):=\int_0^{\cdot}  \phi d^\pi\bar\omega.\qquad  
\ea
 which is in fact an isometry:
\begin{proposition}[Isometry property]
$I_{\bar\omega}\colon \left(\mathbb{V}_a(\bar\omega), \|\cdot\|_{L^2([0,T], a)}\right) \to   \left(\mathcal{H}_a(\bar{\omega}), \|\cdot\|_{\pi}\right)$ is an injective isometry.  \label{prop.quotient}
\end{proposition}

\section{Pathwise nature of the integral}\label{sec:uniqueness}

The definition \eqref{eq.pathwiseintegral} of the  integral $\int \nabla_\omega F(.,\omega) d^\pi\omega$  involves values of the integrand $\nabla_\omega F$ along piecewise constant approximations of $\omega$.
To be able to intepret $\int \nabla_\omega F(.,\omega) d^\pi\omega$  as a {\it pathwise } integral, we must show that it only depends on the values $\nabla_\omega F(t,\omega), t\in[0,T]$ of the integrand along the path $\omega$. 

Our goal is to provide  conditions  under which the following  property holds:

\begin{property}[Pathwise nature of the integral]
	Let $F_1,\, F_2\in \mathbb{C}^{1,2}_b(\Lambda^d_T)$,  and $\omega\in Q_\pi([0,T],\mathbb{R}) \cap C^{\nu}([0,T])$ for some $0< \nu<1/2$.
If $\nabla_\omega F_1(t, \omega_t)=\nabla_\omega F_2(t, \omega_t),\, \forall t\in[0,T]$ then
	\[
	\int_0^T \nabla_\omega F_1\cdot d^{\pi} \omega=\int_0^T \nabla_\omega F_2\cdot d^{\pi}\omega.
	\]
\end{property}
 To give a precise statement, let us introduce the following assumption:
	\begin{hyp}[Horizontal local Lipschitz property]\label{ass.nablaF}\ \\
		A non-anticipative functional  $F\colon \Lambda^d_T\to V$ with values in a finite dimensional real vector space $V$,
		 satisfies the  horizontal locally Lipschitz property if
	
		$$\forall \omega\in D([0,T],\mathbb{R}^d), \exists\ C>0,\, \eta >0, \forall h\geq 0, \forall t\leq T-h, \forall \omega' \in D([0,T],\mathbb{R}^d), $$
		\[
		\|\omega_{t}-\omega'_{t}\|_\infty<\eta,\, \Rightarrow |F(t+h, \omega'_{t})-F(t, \omega'_{t})|\leq C h.
		\]
\end{hyp}
The above assumption asserts the Lipschitz regularity of $F$ in time as the value of the path $\omega$ is frozen. Note that this property is weaker than horizontal differentiability.

We start with a useful lemma:
\begin{lemma}[Expansion formula for regular functionals]
	\label{theorem.controll} Let $\omega\in C^{\nu}([0,T],\mathbb{R}^d)$ for some ${\nu}\in(1/3, 1/2]$ and $F\in \mathbb{C}^{1,1}_b(\Lambda^d_T,\mathbb{R}^n)$ be a non-anticipative functional such that
	$\nabla_{\omega} F\in \mathbb{C}^{1,1}_b(\Lambda^d_T, \mathbb{R}^{n\times d}), \nabla^2_{\omega} F\in \in \mathbb{C}^{1,1}_b(\Lambda^d_T, \mathbb{R}^{n\times d\times d}),$  $ F, \mathcal{D}F, \nabla_\omega^3 F\in Lip(\Lambda^d_T, \|\cdot\|_{\infty})$ and $\nabla_\omega^3 F$  horizontally locally Lipschitz.  Then
	
	\[
	\begin{split}
	F(s,\omega_s)-F(t,\omega_t)=\nabla_{\omega}F(t,\omega_t)(\omega(s)-\omega(t))+\int_t^s\mathcal{D}F(u,\omega_u)du\\
	+\frac 1 2 \langle\nabla_{\omega}^2F(t,\omega_t), (\omega(s)-\omega(t))\, ^t(\omega(s)-\omega(t))\rangle+ O (|s-t|^{3\nu^2+\nu}),
	\end{split}\]
	as $|s-t|\to 0$, uniformly in $t,s\in[0,T].$
\end{lemma}

\begin{proof}
	For the result we will need to further expand the representation \eqref{eq.Rfyformula}. First, note that the second term in the third line of \eqref{eq.Rfy2} can be written as
	\begin{equation}\label{eq.Rfy3}
	\begin{split}
	\int_t^s\int_t^u \partial^2_{ij} F(r,\lambda_r)d\lambda^j(r)d\lambda^i(u)=\partial^2_{ij} F(t,\lambda_t)\int_t^s\int_t^u d\lambda^j(r)d\lambda^i(u)\\+
	\int_t^s\int_t^u \left(\partial^2_{ij} F(r,\lambda_r)-\partial^2_{ij} F(t,\lambda_t)\right)d\lambda^j(r)d\lambda^i(u)
	\end{split}
	\end{equation}
	Using the Lemma (*)  for $G=\partial^2_{ij}F$ and then Fubini's theorem, we get
	
	\begin{equation}\label{eq.Rfy4}
	\begin{split}
	\int_t^s\int_t^u \left(\partial^2_{ij} F(r,\lambda_r)-\partial^2_{ij} F(t,\lambda_t)\right)d\lambda^j(r)d\lambda^i(u)\\=
	\int_t^s\int_t^u \int_t^r \mathcal{D}\partial^2_{ij} F(\tau,\lambda_{\tau})d\tau d\lambda^j(r)d\lambda^i(u)\\+\int_t^s\int_t^u \int_t^r \partial^3_{ijk} F(\tau,\lambda_{\tau})d\lambda_k(\tau)d\lambda^j(r)d\lambda^i(u)\\
	=
	\int_t^s\mathcal{D}\partial^2_{ij} F(\tau,\lambda_{\tau})\Lambda^{ji}_{\tau,s}d\tau+
	\int_t^s\partial^3_{ijk} F(\tau,\lambda_{\tau})\Lambda^{ji}_{\tau,s}\dot{\lambda}^k(\tau)d\tau,
	\end{split}
	\end{equation}
	where
	\[
	\Lambda^{ji}_{a,b}=\int_a^b\int_a^u d\lambda^j(r)d\lambda^i(u)=\int_a^b\left(\lambda^j(u)-\lambda^j(a)\right)\dot{\lambda}^i(u)du.
	\]
	Combining \eqref{eq.Rfy1}, \eqref{eq.Rfy2}, \eqref{eq.Rfy3} and \eqref{eq.Rfy4}  we obtain	
	\begin{equation}\label{eq.Rfyfinal}
	\begin{split}
	R^{F}_{t,s}(\lambda)=\int_t^s\mathcal{D}F(u,\lambda_u)du+\int_t^s\mathcal{D}\partial_i F(r,\lambda_r)\delta\lambda^i_{r,s}dr+\partial_{ij}^2 F(t,\lambda_t)\Lambda^{ji}_{t,s}\\+
	\int_t^s\mathcal{D}\partial^2_{ij} F(\tau,\lambda_{\tau})\Lambda^{ji}_{\tau,s}d\tau+
	\int_t^s\partial^3_{ijk} F(\tau,\lambda_{\tau})\Lambda^{ji}_{\tau,s}\dot{\lambda}^k(\tau)d\tau.
	\end{split}
	\end{equation}
	
	Finally exploiting the symmetry of second and third order derivatives: $\partial^2_{ij} F=\partial^2_{ji} F$,\,  $\partial^3_{ijk} F=\partial^3_{jik} F$, from \eqref{eq.Rfyfinal} we obtain:
	
	\begin{equation}\label{eq.errorformula}
	\begin{split}
	R^{F}_{t,s}(\lambda)=\int_t^s\mathcal{D}F(u,\lambda_u)du+\int_t^s\mathcal{D}\partial_i F(r,\lambda_r)\delta\lambda^i_{r,s}dr+\partial_{ij}^2 F(t,\lambda_t)\tilde{\Lambda}^{ji}_{t,s}\\+
	\int_t^s\mathcal{D}\partial^2_{ij} F(\tau,\lambda_{\tau})\tilde{\Lambda}^{ji}_{\tau,s}d\tau+
	\int_t^s\partial^3_{ijk} F(\tau,\lambda_{\tau})\tilde{\Lambda}^{ji}_{\tau,s}\dot{\lambda}^k(\tau)d\tau,
	\end{split}
	\end{equation}
	where $\tilde{\Lambda}$ is the symmetric part of matrix $\Lambda$. By integration by parts it is easy to compute
	\[
	\tilde{\Lambda}^{ji}_{a,b}=\frac{	\Lambda^{ji}_{a,b}+	\Lambda^{ij}_{a,b}}{2}=\frac 1 2 (\lambda^j(b)-\lambda^j(a))(\lambda^i(b)-\lambda^i(a)).
	\]
	Therefore,
	\begin{equation}\label{eq.errorformula2}
	\begin{split}
	R^{F}_{t,s}(\lambda)=\int_t^s\mathcal{D}F(u,\lambda_u)du+\int_t^s\mathcal{D}\partial_i F(r,\lambda_r)\delta\lambda^i_{r,s}dr+\frac 1 2 \partial_{ij}^2 F(t,\lambda_t)\delta\lambda^i_{t,s} \delta\lambda^j_{t,s}\\+
	\frac 1 2\int_t^s\mathcal{D}\partial^2_{ij} F(\tau,\lambda_{\tau})\delta\lambda^i_{\tau,s}\delta\lambda^j_{\tau,s}d\tau+\frac 1 2
	\int_t^s\partial^3_{ijk} F(\tau,\lambda_{\tau})\delta\lambda^i_{\tau,s}\delta\lambda^j_{\tau,s}\dot{\lambda}^k(\tau)d\tau,
	\end{split}
	\end{equation}
	we decompose the last term in the above as follows
	\[
	\begin{split}
	\int_t^s\partial^3_{ijk} F(\tau,\lambda_{\tau})\delta\lambda^i_{\tau,s}\delta\lambda^j_{\tau,s}\dot{\lambda}^k(\tau)d\tau=
	\int_t^s\partial^3_{ijk} F(t,\lambda_t)\delta\lambda^i_{\tau,s}\delta\lambda^j_{\tau,s}\dot{\lambda}^k(\tau)d\tau\\+ 	\int_t^s \left(\partial^3_{ijk}F(\tau,\lambda_{\tau})-\partial^3_{ijk}F(t,\lambda_t)\right)\delta\lambda^i_{\tau,s}\delta\lambda^j_{\tau,s}\dot{\lambda}^k(\tau)d\tau
	\end{split}	
	\]
	Exploiting the cyclical symmetry of tensor $\partial^3_{ijk} F$ in the indices $i, j, k$:
	\[
	\begin{split}
	\partial^3_{ijk} F(t, \lambda_t)\delta\lambda^i_{\tau,s}\delta\lambda^j_{\tau,s}\dot{\lambda}^k(\tau)=-
	\partial^3_{ijk} F(t, \lambda_t)\delta\lambda^i_{\tau,s}\delta\lambda^j_{\tau,s} \frac{d}{d\tau}\lambda^k_{\tau, s}\\=
	\frac 1 3 \partial^3_{ijk} F(t, \lambda_t)[\delta\lambda^i_{\tau,s}\delta\lambda^j_{\tau,s} \frac{d}{d\tau}\lambda^k_{\tau, s}+\delta\lambda^j_{\tau,s}\delta\lambda^k_{\tau,s} \frac{d}{d\tau}\lambda^i_{\tau, s}+\delta\lambda^k_{\tau,s}\delta\lambda^i_{\tau,s} \frac{d}{d\tau}\lambda^j_{\tau, s}]\\=-\frac 1 3 \partial^3_{ijk} F(t, \lambda_t) \frac{d}{d\tau}\left(\delta\lambda^i_{\tau,s}\delta\lambda^j_{\tau,s} \lambda^k_{\tau, s}\right)
	\end{split}
	\]
	hence 
	\[
	\int_t^s\partial^3_{ijk} F(t,\lambda_t)\delta\lambda^i_{\tau,s}\delta\lambda^j_{\tau,s}\dot{\lambda}^k(\tau)d\tau=
	\frac 1 3 \partial^3_{ijk} F(t, \lambda_t) \delta\lambda^i_{t,s}\delta\lambda^j_{t,s} \lambda^k_{t, s}.
	\]
	Thanks to this \eqref{eq.errorformula2} can be written as
	\begin{equation}\label{eq.errorformula3}
	\begin{split}
	R^{F}_{t,s}(\lambda)=\int_t^s\mathcal{D}F(u,\lambda_u)du+\int_t^s\mathcal{D}\partial_i F(r,\lambda_r)\delta\lambda^i_{r,s}dr+\frac 1 2 \partial_{ij}^2 F(t,\lambda_t)\delta\lambda^i_{t,s} \lambda^j_{t,s}\\+
	\frac 1 2\int_t^s\mathcal{D}\partial^2_{ij} F(\tau,\lambda_{\tau})\delta\lambda^i_{\tau,s}\delta\lambda^j_{\tau,s}d\tau+
	\frac 1 6 \partial^3_{ijk} F(t, \lambda_t) \delta\lambda^i_{t,s}\delta\lambda^j_{t,s} \lambda^k_{t, s}\\
	+\frac 1 2\int_t^s \left(\partial^3_{ijk}F(\tau,\lambda_{\tau})-\partial^3_{ijk}F(t,\lambda_t)\right)\delta\lambda^i_{\tau,s}\delta\lambda^j_{\tau,s}\dot{\lambda}^k(\tau)d\tau.
	\end{split}
	\end{equation}
	Now, as in the proof of part a), we take $\lambda=\omega^N$ and estimate $R^{F}_{t,s}(\omega^N)$.
	To estimate the last term in \eqref{eq.errorformula3} we use that $\nabla^3 F$ is in $Lip([0,T], \|\cdot \|_{\infty})$ and is  locally horizontally Lipschitz:
	\[
	\begin{split}
	|\partial^3_{ijk}F(\tau,\omega^N_{\tau})-\partial^3_{ijk}F(t,\omega^N_t)|\leq |\partial^3_{ijk}F(\tau,\omega^N_{\tau})-\partial^3_{ijk}F(\tau,\omega^N_{t,\tau-t})|\\
	+|\partial^3_{ijk}F(\tau,\omega^N_{t,\tau-t})-\partial^3_{ijk}F(t,\omega^N_t)|\leq C_F \|\omega^N_{t,\tau-t}-\omega^N_{\tau}\|_{\infty}+C_F|\tau-t|\\
	\leq C_{F, \|\omega\|_{\nu} } |\tau-t|^{\nu}+C_F|\tau-t|\leq C_{F,\|\omega\|_{\nu}, T} |s-t|^{\nu}.
	\end{split}
	\] 
	Plugging this estimate into the formula \eqref{eq.errorformula3} and using the local boundedness of the derivatives $\mathcal{D}F, \mathcal{D}\partial_i F,$ $ \partial_{ij}^2 F, \mathcal{D} \partial_{ij}^2 F, \partial^3_{ijk}F$ and the estimates \eqref{eq.omgaNbounds} we obtain
	\begin{equation}\begin{split}
	\left|R^{F}_{t,s}(\omega^N) -\int_t^s\mathcal{D}F(u,\omega^N_u)du- \frac 1 2 \partial_{ij}^2 F(t,\omega_t)\delta\omega^i_{t, s} \delta\omega^j_{t, s}\right|\\
	\leq C_{F, \|\omega\|_{\nu}} |s-t|^{1+\nu}+C_{F, \|\omega\|_{\nu}} |s-t|^{1+2\nu}\\
	+C_{F, \|\omega\|_{\nu}} |s-t|^{3\nu}+ C_{F, \|\omega\|_{\nu}} N^{1-\nu}|s-t|^{4\nu}\\
	\leq C_{F, \|\omega\|_{\nu},T} |s-t|^{3\nu}+C_{F, \|\omega\|_{\nu}} N^{1-\nu}|s-t|^{4\nu}, 
	\end{split}
	\end{equation}
	where we have used also that $\omega^N_t=\omega_t,\, \omega^N(s)=\omega(s)$.
	Recalling \eqref{eq.Rfyapprox} we conclude by triangle inequality
	\begin{equation}\begin{split}
	\left|R^{F}_{t,s}(\omega) -\int_t^s\mathcal{D}F(u,\omega_u)du- \frac 1 2 \partial_{ij}^2 F(t,\omega_t)\delta\omega^i_{t, s} \delta\omega^j_{t, s}\right|\\
	\leq  \int_t^s\left|\mathcal{D}F(u,\omega^N_u)-\mathcal{D}F(u,\omega_u)\right|du+C_FN^{-\nu}|s-t|^{\nu}\\+C_{F, \|\omega\|_{\nu},T} |s-t|^{3\nu}
	+C_{F, \|\omega\|_{\nu}} N^{1-\nu}|s-t|^{4\nu}, 
	\end{split}
	\end{equation}
	The above inequality holds for any $N>1$,thus we can take $N\approx|s-t|^{-3\nu}$ to get
	\[
	\begin{split}
	\left|R^{F}_{t,s}(\omega) -\int_t^s\mathcal{D}F(u,\omega_u)du- \frac 1 2 \partial_{ij}^2 F(t,\omega_t)\delta\omega^i_{t, s} \delta\omega^j_{t, s}\right|\\
	\leq  \int_t^s\left|\mathcal{D}F(u,\omega^N_u)-\mathcal{D}F(u,\omega_u)\right|du
	+C_{F, \|\omega\|_{\nu}} |s-t|^{\nu +3\nu^2}\\
	\leq C_F\|\omega\|_{\nu}|s-t|^{1+\nu}+C_{F, \|\omega\|_{\nu}} |s-t|^{\nu +3\nu^2}, 
	\end{split}\]
	where in the last inequality we used that $\mathcal{D}F$ is Lipschitz continuous.  Hence the result.

\end{proof}

We are now ready to prove the main result of this section, which gives sufficient conditions on the functional $F$ under which the pathwise integral $\int\nabla_\omega F(t,\omega) d^\pi\omega$
is a (pathwise) limit of Riemann sums computed along $\omega$, rather than approximations of $\omega$:

\begin{theorem}\label{theorem.isomuniq}
	Let  $\omega\in Q_\pi([0,T],\mathbb{R}^d)\cap C^{\nu}([0,T],\mathbb{R}^d)$ with $\nu>\frac{\sqrt{13}-1}{6}$ satisfying Assumption \ref{ass.oscill}.
Assume $F \in \mathbb{C}^{1,2}_b(\Lambda^d_T)$ is such that $\nabla_{\omega} F, \nabla^2_{\omega} F,\in \mathbb{C}^{1,1}_b(\Lambda^d_T) $,  $ F, \mathcal{D}F, \nabla^3 F\in Lip(\Lambda^d_T, \|\cdot\|_{\infty})$ and $\nabla_\omega^3 F$ is horizontally locally Lipschitz.  Then the pathwise  integral \eqref{eq.pathwiseintegral} is a limit of non-anticipative Riemann sums:
	\[
	\int_0^T\nabla_{\omega}F(u, \omega) d^{\pi}\omega(u)=\lim_{n\to+\infty} \sum_{[t,s]\in\pi^n}\nabla_{\omega}F(t, \omega)(\omega(s)-\omega(t).
	\]
	In particular, its value only depends on $F(.,\omega)$ and
	\[
	\nabla_{\omega}F(\cdot, \omega) = 0 \implies    \int_0^T\nabla_{\omega}F(u, \omega) d^{\pi}\omega(u)=0.
	\]
\end{theorem}

\begin{proof}
	Let us denote
	\[
	\begin{split}
	A^n_i:=F(t^n_{i+1},\omega_{t^n_{i+1}})-F(t^n_i, \omega_{t^n_i})-\nabla_{\omega}F(t,\omega_{t^n_i})(\omega(t^n_{i+1})-\omega(t^n_i))\\
	-\int_{t^n_i}^{t^n_{i+1}}\mathcal{D}F(u, \omega)du
	-\frac 1 2 \langle \nabla^2_{\omega}F(t,\omega_{t^n_i}),  \omega_{t^n_i, t^n_{i+1}}\, ^t\omega_{t^n_i, t^n_{i+1}}\rangle
	\end{split}
	\]
	From Assumption \eqref{ass.oscill} $\max_i |A^n_i|\to_{n\to\infty}0$ and  by part b) of Theorem \eqref{theorem.controll} $|A^n_i|\leq C |t^n_{i+1}-t^n_i|^{3\nu^2+\nu}$. From the condition on $\nu,$ $3\nu^2+\nu>1$, therefore
	\[
	\sum_i\left|A^n_i \right|\leq C (\max_i |A^n_i|)^{1-\frac{1}{3\nu^2+\nu}}\sum_i|t^n_{i+1}-t^n_i|\leq C T(\max_i |A^n_i|)^{1-\frac{1}{3\nu^2+\nu}}\to 0,
	\]
	which implies $\sum_i A^n_i\to 0$. Hence the limit of Riemann sums exist and is equal to
	\[
	\begin{split}
	\lim_{n\to\infty} \sum _{i}\nabla_{\omega}F(t,\omega_{t^n_i})(\omega(t^n_{i+1})-\omega(t^n_i))
	=F(T,\omega)-F(0,\omega)-\int_0^T \mathcal{D}F(u, \omega)du\\
	-\frac 1 2\int_0^T  \langle\nabla^2_{\omega}F(t,\omega), d[\omega]^{\pi}(t)\rangle= \int_0^T\nabla_{\omega}F(t,\omega) d^\pi\omega(t).
	\end{split}
	\]
	
\end{proof}

\section{Rough-smooth decomposition for functionals of irregular paths}\label{sec:decomposition}

As an application of the previous results, we now derive a decomposition theorem for functionals of irregular paths.
Let $\pi=\{\, \pi^n\,\}$ be a sequence of partitions with $|\pi^n|\to0$. Consider   an irregular path with strictly increasing  quadratic variation along $\pi$:
\ba\bar\omega\in Q_\pi([0,T],\mathbb{R}^d)\cap C^{1/2-}([0,T],\mathbb{R}^d)\quad {\rm with}\quad \frac{d[\bar\omega]_\pi}{dt}>0\ dt-a.e.\label{eq.irregularity}\ea
 We define the vector space of paths obtained as the image of $\bar\omega$ under regular non-anticipative functionals:
$$\mathcal{U}(\bar\omega):=\left\{\,F(\cdot, \bar\omega)\,\big|\, F \ {\rm satisfies\ the\ assumptions\ of\ Theorem\  \ref{theorem.isomuniq}}\  \right\} \subset Q_\pi([0,T],\mathbb{R}).$$
Denote
$$\mathbb{V}_{a}(\bar\omega):=\left\{\,\nabla_{\omega}F(\cdot, \bar\omega)\,\big|\, F \ {\rm satisfies\ the\ assumptions\ of\ Theorem \ \ref{theorem.isomuniq}}\  \right\} \subset Q_\pi([0,T],\mathbb{R}^d).$$
The following result gives a 'signal plus noise' decomposition for such paths, in the sense of Dirichlet processes \cite{follmer1981b}:
\begin{proposition}[Rough-smooth decomposition of paths] Any path $\omega\in \mathcal{U}(\bar\omega)$ has a unique decomposition 
	\ba \omega(t) &=& \omega(0)+\int_0^t \phi d^\pi\bar\omega + s(t)\quad{\rm where}\ \phi\in \mathbb{V}_a(\bar\omega),\qquad  [s]_\pi=0. \label{eq.pathdecomposition}\ea
	\label{prop.decomposition}
\end{proposition}
The term $\int_0^t \phi .d^\pi\bar\omega$ is the 'rough' component of $\omega$ which inherits the irregularity of $\bar\omega $ while $s(.)$ represents a 'smooth' component with zero quadratic variation. 
\begin{remark}[Pathwise 'Doob-Meyer' decomposition]
This result may be viewed as a pathwise analogue of the well-known decomposition of a continuous semimartingale into a local martingale (a process with stricly increasing quadratic variation
A similar pathwise decomposition result was obtained by Hairer and Pillai \cite{hairerpillai2013} using rough path techniques (see also \cite[Theorem 6.5]{frizhairer} and \cite{friz2013}). Unlike the semimartingale decomposition, the result of Hairer and Pillai \cite{hairerpillai2013} involves H\"older-type regularity assumptions on the components, as well as a uniform H\"older roughness condition on the path. 

Our setting is closer to the original semimartingale decomposition in that the components are distinguished based on quadratic variation and the irregularity condition \eqref{eq.irregularity} on $\omega$ is also expressed in terms of quadratic variation.
\end{remark}

\begin{proof} Let $\omega\in \mathcal{U}(\bar\omega).$ Then there exists $F\in \mathcal{R}(\Lambda^d_T)$ verifying the assumptions of Theorem \ref{theorem.isomuniq}, with $\omega(t)=F(t,\bar\omega)$. The  functional change of variable formula applied to $F$ then yields the decomposition with $\phi=\nabla_\omega F(.,\bar\omega)\in \mathbb{V}_a(\bar\omega)$ and $s(t)=F(0,\omega_0)+\int_0^t du ({\mathcal D}F(u,\omega)+\frac{1}{2}<\nabla^2_\omega F(u,\omega),d[\omega]>)$. The continuity of ths integrand implies that $s$ has finite variation so $[s]_\pi=0$. 

Consider now two different decompositions
\[
\omega(t)-\omega(0) = \int_0^t \phi_1d^\pi\bar\omega + s_1(t)= \int_0^t \phi_2d^\pi\bar\omega + s_2(t).
\]
so $s_1-s_2=\int_0 (\phi_1-\phi_2) d^\pi\bar\omega.$
Since $\phi_1-\phi_2\in \mathbb{V}_a(\bar\omega)$ this representation shows that $s_1-s_2\in Q_\pi([0,T],\mathbb{R})$ so the pathwise quadratic covariation $[s_1,s_2]_\pi$ is well defined and we can apply the polarization formula to conclude that  
 $[s_1-s_2]_\pi=0$. Applying Proposition \ref{0quadratic} to $s_1-s_2=\int_0 (\phi_1-\phi_2) d^\pi\bar\omega.$ then implies that that $\phi_1=\phi_2$. 

Finally, applying Theorem \ref{theorem.isomuniq} to $\int_0^t (\phi_1-\phi_2)d^\pi\bar\omega$ then shows that $s_1=s_2$ which  yields uniqueness of the decomposition \eqref{eq.pathdecomposition}. 
\end{proof}
\begin{remark} As noted by Schied \cite{schied2016}, $Q_\pi([0,T],\mathbb{R})$ is not a vector space and, given two paths $(\omega_1,\omega_2)\in Q_\pi([0,T],\mathbb{R})$ the quadratic covariation along $\pi$ cannot be defined in general.
By contrast, the space
$\mathcal{U}(\bar\omega)$ introduced above {\it is} a vector space of paths with finite  quadratic variation along $\pi$.  Moreover,
 for any pair of elements $(\omega_1,\omega_2)\in \mathcal{U}(\bar\omega)^2,$ the quadratic covariation along $\pi$ is well defined; if $\omega_i=\int_0 \phi_i.d^\pi\overline{\omega}+ s_i$ is the rough-smooth decomposition of $\omega_i$ the quadratic covariation is given by
\ba [\omega_1,\omega_2]_\pi(t)
=\int_0^t <\phi_1^t\phi_2,d[\overline{\omega}]>.\ea
\end{remark}

\end{document}